# Analysis of different Mathematical Models for Different Case Studies Using Statistical Fitting

Hamidreza Moradi[1]*

Hamideh Hossei[2]

Erfan Kefayat[3]

[1]Department of Mechanical Engineering and Engineering Science, The University of North Carolina at Charlotte,

[2]Department of Infrastructure and Environmental System, The University of North Carolina at Charlotte,

[3]Department of Public Policy, The University of North Carolina at Charlotte,

Charlotte, North Carolina, USA

*Corresponding Author

**Abstract**

Curve fitting is a fundamental technique in engineering and scientific research, serving as a critical tool for extracting insights from data. This study explores the application of various statistical equations to estimate outcomes in three distinct case studies: population dynamics, temperature variations within buildings, and market equilibrium prices. The efficacy of each fitting is evaluated through rigorous error criteria, including Sum of Squares Error (SSE), R-squared ($R^2$), Degrees of Freedom Error (DFE), Adjusted R-squared (Adj. $R^2$), and Root Mean Square Error (RMSE). Our findings reveal that the selection of mathematical functions and the order of equations are contingent upon the specific nature of the model being analyzed. In the first case study concerning population dynamics, a fractional exponential function emerges as the optimal equation. Conversely, the second and third case studies, which focus on temperature fluctuations in buildings and market price equilibrium, respectively, find the best fit with a sinusoidal function employing three terms. Additionally, we compare the fitted curves and empirical equations with experimental data from the literature to comprehensively assess their predictive accuracy. By establishing the most suitable mathematical models for diverse case studies, researchers and practitioners can enhance their ability to make accurate predictions and informed decisions.

## 1- Introduction

The Ordinary Differential Equation (ODE) method has many applications in various fields of science and engineering, including biology, economics, and physics. Many models have been developed to simulate and predict the growth of biological systems, including animals and humans. The models address different aspects of population dynamics, with some models employing discrete modeling for certain instances and mainly continuous modeling for larger populations (Tsoularis & Wallace, 2002). (Wei & Zhang, 2019) introduced the concept of a stable population with a saturation level known as the carrying capacity, denoted by K, which acts as an upper limit for growth. As a response to this limitation, he developed a logistic growth equation, based on the exponential model(Zabadi, Assaf, & Kanan, 2017).

Engineering projects on material investigation, electronics engineers, mechanical engineers, and chemical engineering problems commonly encounter complex relationships between variables that cannot be captured by simple linear models; thus, Numerical models can take the lead in this case (Bagheri, Siavashi, & Mousavi, 2023). (Moradi, Beh Aein, & Youssef, 2021) uses Response surface method along with fractional factorial design (FFD) to control 12 variables in an optimization problem. In another work, (Noorbakhsh & Moradi, 2020) used Taguchi method to control different parameters in manufacturing processes. However, the solution costs of both works are not balanced with engineering problems considered in this study. (Moradi, Rafi, & Muslim, 2021) uses a novel numerical method to solve the Paris law ODE for the crack growth behavior under dynamic loading. Similarly, this method needs a huge amount of calculation which is not a suitable and optimized method for the cases introduced here. As a result of curve fitting, they were able to uncover complex patterns, enabling accurate predictions, optimizations, and decisions to be made (Mosallaei, Afrouzmehr, & Abshari, 2024). Additionally, curve fitting facilitates parameter estimation and model validation. The alignment of experimental data with theoretical models allows engineers to refine their theories, validate assumptions, and refine their designs.

Numerous studies have used different versions of the exponential model to predict population changes in poultry egg production (Yang, Wu, & McMillan, 1989). The Malthusian model is used to estimate population growth based on these two points (Weir, 1991). The logistic model can be beneficial for more precise projections, since it incorporates three different types of data (Tu, 1996). In addition to biology, numerous other fields have found application to the logistic model. (Fisher & Pry, 1971) effectively utilized this model to describe the market penetration of new products and technologies. (Marchetti, 1977) employed the logistic model to summarize worldwide energy usage and source substitution. Furthermore, (Herman & Montroll, 1972) demonstrated that even fundamental evolutionary processes, like the industrial revolution, can be modeled using logistic dynamics. In addition, they observed a decline in the percentage of the labor force engaged in agriculture as well as an increase in the percentage engaged in industry.

As another case study, the thermal comfort of educational and residential buildings has been extensively studied in previous studies, as these are the primary indoor environments in which people spend a substantial amount of time. However, the focus of those studies has primarily been on specific types of educational facilities, such as kindergartens (Heracleous & Michael, 2019), pre-schools (Lu, Deng, Li, Sundell, & Norbäck, 2016), day care centers (Shriram, Ramamurthy, & Ramakrishnan, 2019), nurseries (Branco, Alvim-Ferraz, Martins, & Sousa, 2015), primary schools (Seppänen, Fisk, & Mendell, 1999), and elementary schools (Duarte, Gomes, & Rodrigues, 2018). Limited attention has been given to exploring the thermal comfort in university buildings (Wargocki, Wyon, Baik, Clausen, & Fanger, 1999). Additionally, many of these studies have overlooked the distinction between thermal comfort experienced during occupancy and non-occupancy periods, as they solely monitored the buildings during working (Quang, He, Knibbs, De Dear, & Morawska, 2014). To the best of the authors' knowledge, none of these studies have investigated the variations in thermal comfort parameters between different types of terms and sources in buildings. As a second case study, we intend to examine how effectively a numerical model can predict changes in a building's internal temperature.

The equilibrium problem in financial economics can be divided into two stages. First, the objective is to determine the optimal portfolio for each agent, considering the current prices of basic securities. Second, the equilibrium price of each basic security is examined by ensuring that total supply equals total demand. The study of determining the current price of an asset, given its future payoff, is known as asset pricing. Notably, the payoff of any derivative security is contingent upon the underlying securities. In the third case study, our focus is on the equilibrium price of the market, aiming to establish more accurate relationships that represent the reference data under consideration.

Testing specific hypotheses based on the parameter values of the best fitting function not only provides insights into underlying processes but also uncovers new knowledge. A diverse array of probabilistic functions can be utilized to describe a set of response times. While some functions have a solid theoretical justification, others present parameter values that are more easily interpretable(Archontoulis & Miguez, 2015; Ferguson, 2014). (Huang, Lu, & Sellers, 2007) suggested a framework suitable for large scale systems in which the agents are weakly linked with each other only through their cost function. This approach finds applications in diverse scenarios, including dynamic economic models involving agents linked via a market and wireless power control. Due to the weak coupling between dynamics and costs, considerable efforts have been dedicated to developing numerical approximations. (Laitner & Stolyarov, 2003) proposed a framework that enables integration of financial market values into economic growth study. The researchers assume that frontier technology evolves over time in an exogenous and uneven manner. Each unit of capital, whether tangible or intangible, embodies a specific technology. Moreover, firms in this framework produce utilizing physical capital, labor, and applied knowledge. In this model, transforming inventions emerge at randomly spaced intervals, and the empirical section assumes that historians and other commentators can identify the key arrival dates.

An alternative equilibrium approach is developed by (Blume, Easley, & O'hara, 1994) to study the behavior of security markets. The model assumes that a fundamental aspect is unknown to all traders, and they receive signals informing them about the asset's fundamental. However, in this model, aggregate supply remains fixed. The noise in the model originates from the quality of information, specifically the precision of the signal distribution. Prices alone cannot provide complete information on both the magnitude of the signals and their precision. Their research demonstrates that volume carries information about the quality of traders' information that cannot be inferred from price statistics. Additionally, sequences of volume and prices can provide informative insights. Traders who utilize information from the market statistics outperform those who do not.

The evaluation of statistical models has been a central concern throughout various fields, leading to the application of diverse metrics, including Sum of Squares Error (SSE), R-squared ($R^2$), Degrees of Freedom Error (DFE), Adjusted R-squared (Adj. $R^2$), and Root Mean Square Error (RMSE). The mentioned metrics offer insights into model performance and fit, with each serving distinct purposes.

(Buse, 1973) introduced extensions of $R^2$ for models with heteroscedastic errors and known variance. This marked a significant step towards accommodating variations in error structures beyond the classical linear regression framework. Subsequently, some notable exceptions emerged, such as the application of $R^2$ to logit and probit models by (Windmeijer, 1995), and tobit models, which were surveyed by (Veall & Zimmermann, 1996). These extensions helped broaden the use of $R^2$ to models with different underlying assumptions and distributions. the application of $R^2$ has evolved from its original linear regression context to encompass models with varying error structures and distributions. While its usage remains sparse in certain fields, its adaptability to different scenarios underscores its enduring value as a metric for assessing model fit, albeit within specific modeling paradigms.

## 2- Methodology

This study delves into the empirical relationships concerning three distinct types of change: population dynamics, temperature change, and market equilibrium price change. The objective is to derive empirical relationships for each of these domains. To achieve this, we meticulously conducted data collection to ensure the reliability and comprehensiveness of our dataset.

### 2-1 Case Studies

In this section, we conduct in-depth case studies to explore and model three distinct types of change: population dynamics, temperature change in a building, and market equilibrium price change. Each case study is grounded in well-established theories and empirical data, providing a comprehensive understanding of these dynamic phenomena.

### 2-1-1 Population Dynamics

Let p(t) represent the population at time t. Although the population is always an integer, it is usually large enough that assuming p(t) is continuous introduces little error. The next step is to determine the growth rate (input) and the death rate (output) for the population.

### Malthusian Model

One of the first models proposed to predict population change is as follows:

$$\frac{dp}{dt} = kp \qquad p(0) = p_0 \tag{1}$$

$$p(t) = p_0 e^{kt} \tag{2}$$

Assuming that t = 0 corresponds to the year 1900, the formula above yields the following result:

$$p(t) = 76.09 e^{kt} \tag{3}$$

In this case, p(t) represents the population in millions. It is possible to fit the model to the data for a specific year (t = 100 years) to obtain a value for k. Finally, we have:

$$p(t) = 76.09 e^{(0.00065528)t} \tag{4}$$

Until about 1900, Malthusian models were reasonably consistent with census data. The predicted population after 1900 is large to be accounted for by the Malthusian model.

### Logistic Model

Moving beyond the Malthusian model, we delve into the more sophisticated logistic model. This model accounts for the limiting factors that eventually slow population growth. As a result, the following equations can be derived:

$$\frac{dp}{dt} = -Ap(p - p_1) \qquad p(0) = p_0 \tag{5}$$

$$p(t) = \frac{p_0 p_1}{p_0 + (p_1 - p_0)e^{-Ap_1 t}} \tag{6}$$

The equation above introduces parameters A, p1, and p0, which we estimate based on empirical data. By fitting the logistic model, we gain a more nuanced understanding of how populations approach equilibrium. It indicates that p₀ = 76.09 when t = 0, which corresponds to the year 1900. Next, we must determine the parameters A, p₁ in the equation above. In order to fit the function with all available data, we obtained constants A = -4.382e-07, p₁ = -2.921e+04, and P₀ = 76.09. The results indicate that:

$$p(t) = \frac{-2.2226e + 06}{76.09 + (-2.9286e + 04)e^{-(0.0128)t}} \tag{7}$$

To enhance our modeling precision, we consider a variety of fitted equations to comprehensively describe population changes. This multi-faceted approach ensures that our analysis is robust and reflective of real-world population dynamics.

**Table 1: Fitted equations for population dynamics**

| Nelder 1961 | $\dfrac{A}{(1+e^{-\frac{(\lambda k \times year)}{\theta}})^\theta}$ | (Nelder, 1961) |
|---|---|---|
| Mc Millan 1980 | $A(e^{-k_2 \times year} - e^{-k_1 \times year})$ | (I. A. N. McMillan, 1980) |
| Mc Millan 1970 | $a \times \dfrac{e^{(-x \times year - c \times year + c \times d)}}{\dfrac{(x+c) - e^{(-x \times year)}}{x + c_1}}$ | (I. McMillan et al., 1970) |
| Mc Nally 1971 | $a \times year^b \times e^{(-c \times year)}$ | (McNally, 1971) |
| Yang 1989 | $\dfrac{a \times e^{(-x \times year)}}{1 + e^{(-c \times (year - d))}}$ | (Yang et al., 1989) |

### 2-1-2- Temperature Change in a Building

The purpose of curve fitting is to compare numerical models that depict the 24-hour temperature profile within a building based on various factors, such as the outside temperature, internal heat generation, and heating or cooling provided by the furnace or air conditioner. We adopt a systematic approach by treating the building as a single compartment in our analysis. Compartmental analysis is an effective method of modeling the interior temperature of a building. In this approach, we represent the building as a single compartment, and the temperature inside the building at time t is denoted by T(t). The rate of change in temperature is determined by all the factors that generate or dissipate heat.

The temperature inside a building is influenced by three main factors. Firstly, there is the heat produced by people, lights, and machines within the building, denoted as H(t). The second component is the heating or cooling provided by the furnace or air conditioner, represented by U(t). As a result of this component, the temperature increases or decreases at a particular rate. Generally, the additional heating rate H(t) and the furnace or air conditioner rate U(t) are described in terms of energy per unit time (e.g., British thermal units per hour). We can, however, express H(t) and U(t) in terms of temperature per unit time by multiplying them by the building's heat capacity (in degrees heat change per unit of heat energy). A third factor influencing the temperature inside a building is the outside temperature, denoted by M(t). Experimental evidence has shown that this factor can be modeled using Newton's law of cooling, which states that the rate of change in the temperature T(t) is proportional to the difference between the outside temperature M(t) and the inside temperature T(t). That is, the rate of change in the building temperature due to M(t) is $K[M(t) - T(t)]$ in which the positive constant K depends on the physical properties of the building, such as the number of doors and windows and the type of insulation, but K does not depend on M, T, or t. In other words, we found:

$$\frac{dT}{dt} = K[M(t) - T(t)] + H(t) + U(t) \tag{8}$$

$$\frac{dT}{dt} + KT(t) = Q(T) \tag{9}$$

Where $Q(T) = KM(t) + H(t) + U(t)$. By solving T(t) using the integrating factor $\mu(t) = e^{(\int k \, dt)} = e^{kt}$ and multiplying both sides of 6 by the calculated factor, we have:

$$T(t) = e^{-kt} \int e^{kt} Q(t) dt + Ce^{-kt} \tag{10}$$

Accordingly, the following equation can be derived from it.

$$T(t) = e^{-kt}\left\{\int e^{kt}[KM(t) + H(t) + U(t)]dt + C\right\} \quad (11)$$

If the additional heating rate H(t) equals the constant $H_0$, then the outside temperature M varies as a sine wave over a 24-hour period, with its minimum at midnight and its maximum at noon. That is $M(t) = M_0 - B\cos(\omega t)$

Let B be a positive constant, $M_0$ represent the average outside temperature, and ω be the angular frequency defined as $\omega = \frac{2\pi}{24} = \frac{\pi}{12}$ radians/hr. Furthermore, suppose a simple thermostat is installed, which compares the actual internal temperature with a desired internal temperature TD. If the actual temperature is lower than the desired temperature, the furnace supplies heating; otherwise, it remains off. Additionally, when the actual temperature exceeds the desired temperature, the air conditioner provides cooling; otherwise, it does not function. For simplicity, we disregard any dead zone around the desired temperature, where the temperature difference is insufficient to activate the thermostat.

Assuming that the amount of heating or cooling supplied is directly proportional to the difference in temperature, meaning that $U(t) = K_U[T_D - T(t)]$, where $K_U$ is the (positive) proportionality constant. Now we have:

$$\frac{dT}{dt} = K[M(t) - T(t)] + H_0 + K_U[T_D - T(t)] \quad (12)$$

When the additional heating rate is a constant $H_0$ and the outside temperature M varies as a sine wave over a 24-hr period, the forcing function is:

$$Q(T) = K(M_0 - B\cos\cos(\omega t)) + H_0 + K_U T_D \quad (13)$$

The following results can be obtained by changing some of the parameters:

$$\omega = \frac{2\pi}{24} = \frac{\pi}{12} \qquad k_1 = k + K_U \quad (14)$$

$$B_2 = \frac{K_U T_D + K M_0 + H_0}{k_1} \qquad B_1 = \frac{BK}{k_1} \qquad (15)$$

Hence, the solution to the differential equation 5 is as follows.

$$T(t) = B_2 - B_1 F_1(t) + Ce^{(-k_1 t)} \qquad (16)$$

Where:

$$F_1(t) = \frac{\cos(\omega t) + \left(\frac{\omega}{k_1}\right)\sin(\omega t)}{1 + \left(\frac{\omega}{k_1}\right)^2} \qquad (17)$$

The constant C is chosen so that at time t = 0 the value of the temperature T equals $T_0$. Thus, $C = T_0 - B_2 + B_1 F_1(0)$.

The numerical model elucidates the correlation between the temperature distribution within a building and both external temperature and various internal factors. Throughout a typical day, the temperature inside a building experiences fluctuation due to several internal influences, such as the presence of people, lighting, window opening/closing, machines, furnaces, and air conditioning. On the other hand, the external temperature follows a sine wave pattern over a 24-hour period, reaching its minimum at midnight (t=0) and its maximum at noon (t=12). To maintain the desired internal temperature, a thermostat is installed. Whenever the actual temperature deviates from the desired level, the thermostat triggers either the heater or air conditioning to make appropriate adjustments. By inserting all the terms into the equation 10 we get the temperature change of a building as follows:

$$T(t) = \left(\frac{K_U \times td + km_0 + h_0}{k + K_U}\right)$$

$$- \left(\left(\frac{b \times k}{k + K_U}\right) \times \left(\frac{\cos\cos\left(t \times \frac{\pi}{12}\right) + \left(\frac{\pi}{12 \times (k + K_U)}\right) \times \sin\sin\left(\frac{t \times \pi}{12}\right)}{1 + \left(\frac{\pi}{12 \times (k + K_U)}\right)^2}\right)\right.$$

$$+ \left(t_0 - \left(\frac{K_U \times td + k \times m_0 + h_0}{(k + K_U)}\right)\right.$$

$$\left.\left.+ \left(\left(\frac{b \times k}{k + K_U}\right) \times \left(\frac{1}{\left(1 + \left(\frac{\pi}{12 \times (k + K_U)}\right)\right)^2}\right)\right)\right)^{e^{(-(k+K_U) \times t)}}\right) \quad (18)$$

Constants can be expressed as $b$ = -1.04, $h_0$ = 86.61, $k$ = 1.11, $K_U$ = -0.898, $m_0$ = -9.56, $t_0$ = 32.09, and $td$ = 76.41.

Multiple fitted equations are considered to ensure a robust representation of temperature profiles within buildings, allowing us to assess the influence of external temperature and internal dynamics effectively.

**Table 2: Fitted equations for temperature change**

| | |
|---|---|
| Exp. & Sin | $(a_0 + b_1(a_1 \times t - a_2)) \times e^{(f \times t)}$ |
| Fourier | $a_0 + a_1 \cos\cos(x \times w) + b_1 \times \sin\sin(x \times w) + a_2 \times \cos\cos(2x \times w) + b_2 \sin\sin(2x \times w)$ |
| Gauss | $a_1 \times e^{\left(-\left(\frac{(x-b_1)}{c_1}\right)^2 + a_2 \times e^{\left(-\left(\frac{(x-b_2)}{c_2}\right)^2\right)}\right)}$ |
| Exp. | $a \times e^{(b \times x)} + c \times e^{(d \times x)}$ |
| Sin | $a_1 \sin\sin(b_1 x + c_1) + a_2 \sin\sin(b_2 x + c_2) + a_3 \sin\sin(b_3 x + c_3)$ |

**Market Equilibrium Price**

Market equilibrium is defined as a situation where the quantity demanded of a commodity is just equal to the quantity supplied, at a specific price. This analysis employs a linear model to examine how market prices evolve in response to changing conditions. The linear model of the market equilibrium assumes that the demand and supply functions have the form $q_d = d_0 - d_1 p$ and $q_s = -s_0 + s_1 p$, respectively, where p is the market price of the product, $q_d$ is the associated quantity demanded, $q_s$ is the associated quantity supplied, and $d_0$, $d_1$, $s_0$, and $s_1$ are all positive constants. The functional forms ensure that the "laws" of downward sloping demand and upward sloping supply are being satisfied. It is easy to show that the equilibrium price is $\dot{p} = \frac{(d_0+s_0)}{(d_1+s_1)}$.

Economists typically assume that markets are in equilibrium and justify this assumption with the help of stability arguments. For example, consider the simple price adjustment equation:

$$\frac{dp}{dt} = \lambda(q_d - q_s) \tag{19}$$

where $\lambda > 0$ is a constant indicating the speed of adjustments. This follows the intuitive requirement that price rises when demand exceeds supply and falls when supply exceeds demand. When $p^*$ is satisfied as going to infinity for every initial price level p(0), the market equilibrium is said to be globally stable. Therefore, the solution would be as follows:

$$p(t) = [p(0) - \dot{p}]e^{ct} + \overset{.}{p}, where\ c = -\lambda(d_1 + s_1) \tag{20}$$

Now consider a model that considers the expectations of agents. Based on the assumption that market demand and supply functions over time t >= 0 are as follows:

$$q_d(t) = d_0 - d_1 p(t) + d_2 \dot{p}(t)\ and\ q_s(t) = -s_0 - s_1 p(t) - s_2 \dot{p}(t) \tag{21}$$

respectively, where p(t) is the market price of the product, $q_d(t)$ is the associated quantity demanded, $q_s(t)$ is the associated quantity supplied, and $d_0$, $d_1$, $d_2$, $s_0$, $s_1$, and $s_2$ are all positive constants. If the price p(t) of a product is $5 at t = 0 months and demand and supply functions are modeled as:

$$q_d(t) = 30 - 2p(t) + 4\dot{p}(t)\ and\ q_s(t) = -20 + p(t) - 6\dot{p}(t) \tag{22}$$

Taking the ODE form as a starting point, we can derive the exact solution as follows:

$$p(t) = \frac{\lambda a - \lambda b p}{1 - \lambda c} \tag{23}$$

$$P(t) = D e^{\frac{\lambda b t}{c\lambda - 1}} + \frac{a}{b} \quad \{a = d_0 + s_0 \; b = d_1 + s_1 \; c = d_2 + s_2 \tag{24}$$

In the case of constants, they should be expressed as $D$ = 3.282e-08, $L$ = -7.314, $a$ = 36.07, $b$ = -0.01, $c$ = 0.018.

To improve the accuracy of our modeling, we explore a range of tailored equations that collectively provide a comprehensive portrayal of equilibrium price. This multifaceted methodology guarantees that our analysis is both resilient and representative of actual market equilibrium price change.

**Table 3: Fitted equations for equilibrium price**

| Distr. & Exp. | $(A t^2 + Bt + C)^{(Ft^G)}$ |
|---|---|
| Rat | $\dfrac{(p_1 t)^2 + (p_2 t + p_3)}{t + q_1}$ |
| Gauss | $a_1 \times e^{\left(-\left(\frac{(x-b_1)}{c_1}\right)^2\right)} + a_2 \times e^{\left(-\left(\frac{(x-b_2)}{c_2}\right)^2\right)}$ |
| Exp. | $a \times e^{(b \times x)} + c \times e^{(d \times x)}$ |
| Sin | $a_1 \sin \sin (b_1 x + c_1) + a_2 \sin \sin (b_2 x + c_2) + a_3 \sin \sin (b_3 x + c_3)$ |

In this study to evaluate the results of each fitting process we consider the following metrics SSE, a foundational metric, quantifies the sum of squared differences between observed and predicted values. It measures the overall model error, highlighting how well the model captures the data's variability. Lower SSE values signify a better fit.

$R^2$ initially developed for linear regression, assesses the proportion of variability in the dependent variable explained by the model. It ranges from 0 to 1, where higher values imply a greater portion of variability captured. However, $R^2$ is primarily effective in linear contexts with homoscedastic errors, limiting its broader applicability.

DFE represents the degrees of freedom associated with the model's residuals, indicating the extent to which the data's complexity is effectively captured by the model. It aids in assessing the model's flexibility to accommodate variations without overfitting.

Adjusted $R^2$, an enhancement of $R^2$, addresses the potential issue of overfitting by accounting for the number of predictors. It penalizes models with excessive predictors that might not truly contribute to explaining variability, providing a more balanced view of model quality.

RMSE gauges the average magnitude of prediction errors. It quantifies the square root of the average squared differences between observed and predicted values. Lower RMSE values signify better prediction accuracy.

### 3- Results and discussion

In the context of our case studies, we have applied the corresponding curves to estimate the outcomes, whether in their experimental or empirical forms. It's important to note that the choice of evaluation metrics holds varying relevance and effectiveness depending on the specific modeling context.

### 3-1- Population

To delve deeper into our analysis of population dynamics, we obtained the coefficients for each equation through rigorous curve fitting. These coefficients serve as critical parameters in understanding and modeling the intricacies of population change.

**Table 4: Constants for the fitted equations of population dynamics**

| Nelder 1961 | A = 391.8 |
| --- | --- |
|  | k = 2.391 |
|  | $\lambda$ = 581.6 |
|  | $\theta$ = 137.8 |
| Mc Millan 1980 | A = 88.42 |
|  | $k_1$ = 1.689 |
|  | $k_2$ = -0.01136 |

| Mc Millan 1970 | a = 1 |
| | c = 0.5 |
| | c1 = 2.571e+04 |
| | d = -80 |
| | x = 3.909e-05 |
| Mc Nally 1971 | a = 65.4 |
| | b = 0.1047 |
| | c = -0.009579 |
| Yang 1989 | a = 87.99 |
| | c = 98 |
| | d = -934 |
| | x = -0.01141 |

Analyzing the outcomes of our population change models, we turn our attention to key evaluation metrics. The employment of metrics such as Sum of Squared Errors (SSE) and Root Mean Square Error (RMSE) proves invaluable in providing a clear understanding of the error magnitudes associated with each model. On the other hand, R-squared ($R^2$) and Adjusted R-squared (Adj. $R^2$) offer insights into the proportion of variability explained by these models.

Upon a comprehensive evaluation of the five equations employed, it becomes evident that all of them exhibit notably high R-squared values. This observation underscores their effectiveness in capturing and representing the underlying data patterns. Equation 5, in particular, stands out with the highest R-squared value of 0.9914, followed closely by Equation 4 with 0.9877, Equation 2 with 0.9844, Equation 1 with 0.9667, and Equation 3 with 0.6940.

Furthermore, the Root Mean Square Error (RMSE) values across all five equations remain consistently low. This is a strong indicator of the accuracy of their predictions. Equation 5 achieves the lowest RMSE value at 7.4113, followed by Equation 4 at 8.8359, Equation 2 at 9.9787, Equation 1 at 14.6213, and Equation 3 at 44.5107.

In light of these assessments, it is evident that Equation 5 not only exhibits the best overall performance but also strikes a harmonious balance between the complexity of the model and the quality of its fit. This conclusion is reinforced by the Degree of Freedom for Error (DFE), further underlining Equation 5 as the most robust choice for modeling population change.

Table 5: Statistical results for the fitted equations of population dynamics

|  | SSE | R Square | dfe | Adj-R Square | RMSE |
|---|---|---|---|---|---|
| Nelder 1961 | 2.5654e+04 | 0.9667 | 120 | 0.9659 | 14.6213 |
| Mc Millan 1980 | 1.2049e+04 | 0.9844 | 121 | 0.9841 | 9.9787 |
| Mc Millan 1970 | 2.3576e+05 | 0.6940 | 119 | 0.6837 | 44.5107 |
| Mc Nally 1971 | 9.4469e+03 | 0.9877 | 121 | 0.9875 | 8.8359 |
| Yang 1989 | 6.5912e+03 | 0.9914 | 120 | 0.9912 | 7.4113 |

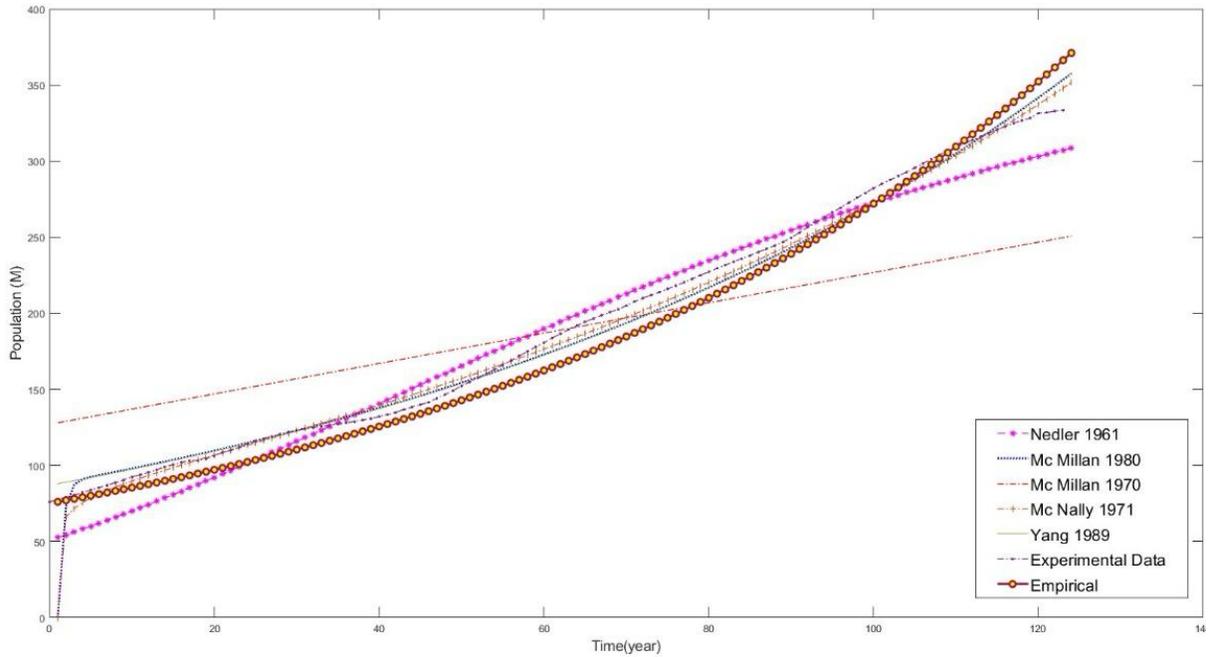

Figure 1. Fitted equations, experimental data, and empirical equation for population dynamics

**3-2- Temperature Change**

For our temperature case study, we delve into the coefficients of each equation, which can be derived through the fitting process as outlined below:

**Table 6: Constants for the fitted equations of temperature change**

| Exp. & Sin | | |
|---|---|---|
| | f | -0.003493 |
| | $a_0$ | 34.6 |
| | $a_1$ | 0.4055 |
| | $a_2$ | 2.636 |
| | $b_1$ | 2.795 |
| Fourier | | |
| | $a_0$ | 33.23 |
| | $a_1$ | -1.456 |
| | $b_1$ | 2.644 |
| | $a_2$ | 0.2962 |
| | $b_2$ | 0.2962 |
| | w | 0.3532 |
| Gauss | | |
| | $a_1$ | 35.94 |
| | $b_1$ | 5.404 |
| | $c_1$ | 16.7 |
| | $a_2$ | 22.41 |
| | $b_2$ | 22.87 |
| | $c_2$ | 6.083 |
| Exp. | | |
| | a | 36.07 |
| | b | -0.01017 |
| | c | 0.01805 |
| | d | 0.2604 |
| Sin | | |
| | $a_1$ | 38.29 |
| | $b_1$ | 0.001099 |
| | $c_1$ | 1.039 |

|  | $a_2$ | 3.117 |
|---|---|---|
|  | $b_2$ | 0.3428 |
|  | $c_2$ | -0.4052 |
|  | $a_3$ | 1.508 |
|  | $b_3$ | 0.7265 |
|  | $c_3$ | 1.175 |

In a similar vein to our population model analysis, we employ Sum of Squared Errors (SSE) and Root Mean Square Error (RMSE) to gauge the error magnitudes associated with each temperature curve. Additionally, we rely on R-squared ($R^2$) values and Adjusted R-squared (Adj. $R^2$) to assess the proportion of variability explained by these temperature models.

Upon careful examination of the five equations utilized in our temperature modeling, it becomes evident that all of them exhibit commendably high R-squared values. This reaffirms their effectiveness in accurately representing the data. Notably, Equation 5 (sin3) stands out with the highest R-squared value of 0.9337, closely trailed by Equation 2 (fourier2) at 0.9279, Equation 3 (gauss2) at 0.7496, Equation 1 (exp and sin) at 0.6374, and Equation 4 (exp2) at 0.3683.

Furthermore, the Root Mean Square Error (RMSE) values for these temperature models consistently demonstrate their ability to make precise predictions. Equation 5 (sin3) leads the pack with the lowest RMSE value of 0.6349, followed by Equation 2 (fourier2) at 1.6387, Equation 3 (gauss2) at 1.1906, Equation 1 (exp and sin) at 1.4494, and Equation 4 (exp2) at 1.8496.

In summary, all five linear regression equations in the table exhibit a strong fit with the temperature data, boasting high R-squared values and low RMSE values. Among these models, Equation 5 (sin3) emerges as the top performer, with the highest R-squared and the lowest RMSE values. Moreover, we consider the Degree of Freedom for Error (DFE), which underscores the equilibrium between model validity and complexity, further reinforcing Equation 5 (sin3) as the preferred choice for modeling temperature changes.

**Table 7: Statistical results for the fitted equations of temperature change**

|  | SSE | R Square | dfe | Adj-R Square | RMSE |
|---|---|---|---|---|---|
| Exp. & Sin | 90.3373 | 0.6374 | 43 | 0.5868 | 1.4494 |
| Fourier | 17.9502 | 0.9279 | 44 | 0.9198 | 1.6387 |
| Gauss | 62.3698 | 0.7496 | 44 | 0.7212 | 1.1906 |
| Exp. | 157.3684 | 0.3683 | 46 | 0.3271 | 1.8496 |
| Sin | 16.5269 | 0.9337 | 41 | 0.9207 | 0.6349 |

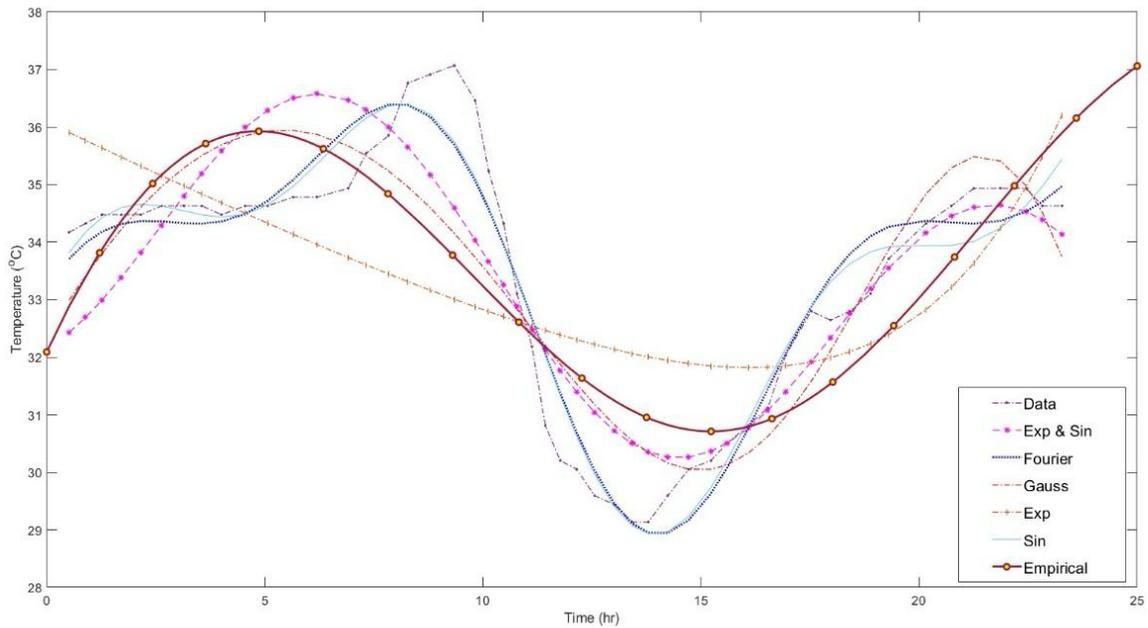

**Figure 2. Fitted equations, experimental data, and empirical equation for temperature change**

### 3-3-Equilibrium Price

Following the curve fitting process, we've derived the coefficients for each equation, which are detailed as follows:

**Table 8: Constants for the fitted equations of equilibrium price**

| Distr. & Exp. | | |
|---|---|---|
| | A | 0.1348 |
| | B | -5.402 |
| | C | 55.84 |
| | F | 0.0109 |
| | G | 1.933 |

| Rat | | |
|---|---|---|
| | p1 | 4.821e+04 |
| | p2 | -9.299e+05 |
| | p3 | 1.074e+07 |
| | q1 | 1.805e+05 |

| Gauss | | |
|---|---|---|
| | a1 | 6.334e+15(-3.089e+19,3.091e+19) |
| | b1 | 129.2 (-1.62e+04, 1.646e+04) |
| | c1 | 19.04 (-1416, 1455) |
| | a2 | 77.81 (-2698, 2854) |
| | b2 | -55.95 (-3741, 3629) |
| | c2 | 79.91 (-2298, 2457) |

| Exp. | | |
|---|---|---|
| | a | 46.84 (39.72, 53.96) |
| | b | -0.01575 (-0.03312, 0.001615) |
| | c | 1.068e-05 (-7.964e-05, 0.000101) |
| | d | 0.7389 (0.3373, 1.14) |

| Sin | | |
|---|---|---|
| | a1 | 495.5 (-1.251e+06, 1.252e+06) |
| | b1 | 0.2139 (-50.18, 50.61) |
| | c1 | -1.561 (-588.5, 585.4) |
| | a2 | 775.4 (-3.472e+05, 3.487e+05) |

|  | b2 | 0.2692 (-102.6, 103.2) |
|  | c2 | 0.9461 (-1160, 1162) |
|  | a3 | 331.3 (-9.081e+05, 9.088e+05) |
|  | b3 | 0.3264 (-49.41, 50.06) |
|  | c3 | 3.452 (-541.7, 548.6) |

Much like our previous models, we apply the same rigorous evaluation criteria to determine the most suitable equation from the table. In this case study, Equation 5, referred to as sin3, emerges as the optimal choice. It boasts the highest R-squared value of 0.9419, signifying an exceptional fit to the data, and it also showcases the lowest RMSE value of 4.1622, indicating remarkable predictive accuracy.

As is evident from our comprehensive literature review, certain metrics, like $R^2$ in Generalized Linear Models (GLMs), hold particular relevance for specific model types. This underscores the importance of meticulous metric selection tailored to the precise analytical objectives and unique characteristics of the dataset under examination.

Table 9: Statistical results for the fitted equations of equilibrium price

|  | SSE | R Square | dfe | Adj-R Square | RMSE |
|---|---|---|---|---|---|
| Distr. & Exp. | 484.5365 | 0.8645 | 16 | 0.8306 | 5.5030 |
| Rat | 1.5589e+03 | 0.5640 | 17 | 0.4871 | 9.5760 |
| Gauss | 620.5019 | 0.8265 | 15 | 0.7686 | 6.4317 |
| Exp. | 594.7089 | 0.8337 | 17 | 0.8043 | 5.9146 |
| Sin | 207.8827 | 0.9419 | 12 | 0.9031 | 4.1622 |

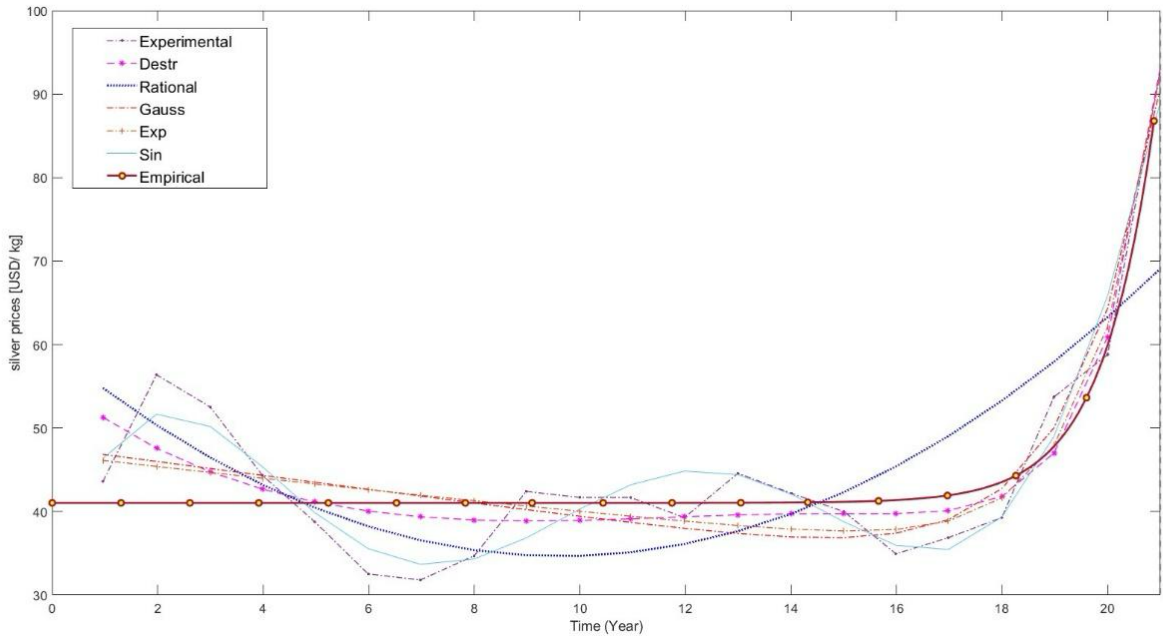

**Figure 3. Fitted equations, experimental data, and empirical equation for equilibrium price**

## 4- Conclusion

In the course of this study, we undertook a comprehensive examination of five distinct equations across three distinct case studies, with the aim of discerning the efficacy of each equation in estimating and predicting results. Our evaluation criteria encompassed five diverse error metrics: Sum of Squared Errors (SSE), Root Mean Square Error (RMSE), R-squared ($R^2$), Adjusted R-squared (Adj. $R^2$), and Degrees of Freedom for Error (DFE).

In the context of our initial case study, involving the logistic model of population dynamics, it became evident that the most fitting curve was a fractional exponential equation. This particular equation, while suggested in prior literature for predicting poultry egg production, exhibited unique constants in our application, emphasizing the importance of tailored model parameterization.

Turning to the subsequent case studies, both the second and third scenarios found optimal fit with a sinusoidal function featuring three terms. Notably, these equations displayed distinct constant values specific to each case, further underscoring the nuanced nature of curve fitting.

To provide a comprehensive comparative analysis, we incorporated experimental data from previous literature into our figures for each case study, facilitating a robust evaluation of empirical solutions against curve fittings.

The implications of this study are multifaceted and extend beyond the realm of curve fitting. It underscores the critical significance of selecting models tailored to the specific attributes of the dataset, advocating for a nuanced approach that acknowledges the unique characteristics of the data under scrutiny. Ultimately, this study contributes to the ongoing refinement of curve fitting methodologies, furnishing valuable insights for diverse fields reliant on data-driven modeling.